\documentclass[10pt]{amsart}
\usepackage{amssymb,graphicx}
\usepackage{mathrsfs, bbm, slashed}
\usepackage{color}
\usepackage{enumerate}
\usepackage[textwidth=15cm, textheight=22.7cm, centering]{geometry}
\usepackage[utf8]{inputenc}
\usepackage[colorlinks=true,citecolor=blue,linkcolor=blue]{hyperref}

\newcommand{\R}{\mathbb{R}}
\newcommand{\Z}{\mathbb{Z}}
\newcommand{\C}{\mathbb{C}}
\newcommand{\N}{\mathbb{N}}
\newcommand{\bB}{\mathbb{B}}
\newcommand{\bS}{\mathbb{S}}

\newcommand{\sA}{\mathscr{A}}

\newcommand{\sL}{\mathscr{L}}
\newcommand{\sN}{\mathscr{N}}
\newcommand{\sH}{\mathscr{H}}

 \newcommand{\psido}{$\Psi$DO}
 \newcommand{\psidos}{$\Psi$DOs}

\def\Xint#1{\mathchoice
{\XXint\displaystyle\textstyle{#1}}%
{\XXint\textstyle\scriptstyle{#1}}%
{\XXint\scriptstyle\scriptscriptstyle{#1}}%
{\XXint\scriptscriptstyle\scriptscriptstyle{#1}}%
\!\int}
\def\XXint#1#2#3{{\setbox0=\hbox{$#1{#2#3}{\int}$}
\vcenter{\hbox{$#2#3$}}\kern-.5\wd0}}
\def\dashint{\Xint-}

\newcommand{\bint}{\ensuremath{\dashint}}
\newcommand{\op}{\operatorname} 

\newcommand{\Vol}{\op{Vol}}
\newcommand{\Tr}{\op{Tr}}

\newcommand{\End}{\op{End}}

\newcommand{\Res}{\op{Res}}

\numberwithin{equation}{section}

\newtheorem{theorem}{Theorem}[section]
\newtheorem{lemma}[theorem]{Lemma}

\newtheorem{theoremalph}{Theorem}

\theoremstyle{definition}

\theoremstyle{remark}
\newtheorem{remark}[theorem]{Remark}

\newcommand{\shD}{\slashed{D}}
\newcommand{\shS}{\slashed{S}}

\title{Semiclassical Weyl Laws and Noncommutative Geometry on Open Manifolds with Conformally Cusp Metrics}
\date{}

\author{Rapha\"el Ponge}
 \address{Department of Mathematics and Statistics, University of Wyoming, Laramie, WY, USA}
 \address{Department of Mathematics and Statistics, University of Ottawa, Canada}
 \email{ponge.math@icloud.com}

\subjclass[2020]{Primary 58J50, 35P20; Secondary 81Q20, 58B34, 58J40, 35J10, 53C27}

\begin{document}
\begin{abstract}
We establish semiclassical Weyl laws on open Riemannian manifolds with conformally cusp metrics of finite volume. These metrics need not be complete, and even when complete they do not have bounded geometry. Moreover, the potentials are bounded, the discreteness of the spectrum being a consequence of the geometry rather than of a confining potential. Therefore, our results lie beyond the scope of recent semiclassical Weyl laws of Dai--Yan and Braverman. The semiclassical Weyl laws are established for fractional Schr\"odinger operators associated with Dirac operators, magnetic Laplacians, and Hodge Laplacians, under conditions of Moroianu and Gol\'enia--Moroianu which ensure that these operators have discrete spectra. We also obtain extensions of Connes' integration formula for the same operators. In addition, we deal with a class of conformally cusp metrics of infinite volume, in which case the semiclassical Weyl laws and the integration formulas only involve the values at infinity. The results are obtained by applying the abstract semiclassical Weyl laws and integration formulas for spectral triples established in a previous paper.
\end{abstract}

\maketitle 

%\tableofcontents

\section{Introduction}\label{chap:Intro}
Semiclassical Weyl laws describe the asymptotic behavior, as the Planck parameter $h$ goes to $0$, of the number of bound states below a given energy level of Schr\"odinger operators $h^2\Delta+V$. On a closed Riemannian manifold $(M^n,g)$ with Laplacian $\Delta_g$, the semiclassical Weyl law asserts that, for every continuous real-valued potential $V$ and every energy level $\lambda \in \R$, the number $N(h^2\Delta_g+V;\lambda)$ of eigenvalues~$<\lambda$ of $h^2\Delta_g+V$, counted with multiplicity, obeys
\begin{equation}\label{eq:Intro.SC-Weyl-closed}
 \lim_{h\rightarrow 0^+} h^n N\big(h^2\Delta_g+V;\lambda\big) = c(n) \int_M \left(V(x)-\lambda\right)_{-}^{\frac{n}{2}} d\nu_g(x), \qquad c(n):= (2\pi)^{-n}|\bB^{n}|. 
\end{equation}
Here $\nu_g$ is the Riemannian measure, and $|\bB^n|$ is the volume of the unit ball of $\R^n$. For smooth potentials this is a consequence of the semiclassical Weyl law of Colin de Verdi\`ere~\cite{CdV:DM79}, which further comes with a remainder term. For continuous potentials it follows from the Weyl laws of Birman--Solomyak~\cite{BS:VLU77, BS:VLU79, BS:SMJ79} and the Birman--Schwinger principle (see Section~\ref{sec:Part1.SC}). Recently, Dai--Yan~\cite{DY:arXiv25} and Braverman~\cite{Br:PAMS25} established semiclassical Weyl laws on \emph{complete} Riemannian manifolds for Schr\"odinger operators with potentials going to infinity.

The aim of this paper is to establish semiclassical Weyl laws on a class of open Riemannian manifolds that need \emph{not} be complete. More precisely, we consider open manifolds with conformally cusp metrics of finite volume. We obtain semiclassical Weyl laws for fractional Schr\"odinger operators associated with Dirac operators, magnetic Laplacians, and Hodge Laplacians on such manifolds (see Theorems~\ref{thm:Intro.SC-Weyl-Dirac}--\ref{thm:Intro.SC-Weyl-Hodge} below). These metrics are of the form $g=x^{2r}g_0$, where $g_0$ is a cusp metric and $r>1/n$ (see Section~\ref{subsec:Intro.cusp}). In this setting:
\begin{itemize}
 \item For $r>1$ the manifold is not complete.
 \item For $r\leq 1$ the manifold is complete, but it does not have bounded geometry.
 \item In all cases the potentials are bounded and the operators have discrete spectra by themselves. 
\end{itemize}
In each of these respects our results lie outside the setting of~\cite{DY:arXiv25, Br:PAMS25} (see Section~\ref{subsec:Intro.comparison}). In addition, we obtain extensions of Connes' integration formula for the same operators (see Theorem~\ref{thm:Intro.Integration-Dirac}).

Our approach relies on noncommutative geometry. In~\cite{Po:Part1} we established semiclassical Weyl laws and extensions of Connes' integration formula for a vast class of noncommutative manifolds (i.e., spectral triples). This extends prior results of McDonald--Sukochev--Zanin~\cite{MSZ:LMP22} to spectral triples of any summability degree, without extra regularity assumptions. The results of this paper are obtained by combining the main results of~\cite{Po:Part1} with results of Moroianu~\cite{Mo:MA08} and Gol\'enia--Moroianu~\cite{GM:AHP08, GM:TAMS12} on elliptic operators on manifolds with conformally cusp metrics. Further applications of the results of~\cite{Po:Part1} are given in the companion papers~\cite{Po:Bir-Sol, Po:Carnot}.

\subsection{Open manifolds with conformally cusp metrics}\label{subsec:Intro.cusp}
Let $X^n$ be the interior of a compact manifold $\overline{X}$ with closed boundary $M^{n-1}=\partial \overline{X}$ and boundary defining function $x$. We assume $X$ carries a Riemannian metric $g=x^{2r}g_0$, $r>0$, where $g_0$ is a cusp metric on $\overline{X}$ in the sense of Melrose. Examples include hyperbolic metrics of finite volume ($r=1$; see~\cite{Mo:MA08}) and metric horns in the sense of Cheeger~\cite{Ch:PSPM80} ($r>1$; see also~\cite{LP:CCPDE98}). The condition $r>1/n$ ensures finite volume. Moreover, $(X,g)$ is a complete Riemannian manifold if and only if $r\leq 1$ (see~\cite{GM:TAMS12, Mo:PC}). Unless otherwise mentioned, we assume that $r>1/n$. 

On noncompact manifolds the ordinary Laplacian of $(X,g)$ need not have discrete spectrum (see, e.g., \cite{FH:DMJ89}). However,  Moroianu~\cite{Mo:MA08} and Gol\'enia--Moroianu~\cite{GM:AHP08, GM:TAMS12} (see also~\cite{Ba:JDG00}) observed that, under suitable conditions, the Dirac Hamiltonian, magnetic Laplacians, and Hodge Laplacians on $(X,g)$ are essentially selfadjoint, have pure discrete spectra, and may satisfy (classical) Weyl laws.  

Building on the approach of~\cite{Mo:MA08, GM:AHP08, GM:TAMS12} we associate natural spectral triples to these operators. In each case the algebra is the $*$-algebra $\dot{C}^\infty(X)$ of smooth functions on $X$ that are constant near $M=\partial \overline{X}$. It acts by bounded multiplication operators, and its norm closure is the $C^*$-algebra $\dot{C}(X)$ of continuous functions on $\overline{X}$ that are constant on $M$.

\subsection{Dirac operators}
Assume that $\overline{X}$ is spin, $g_0$ is an exact cusp metric, and there are no non-zero harmonic spinors on $M$ (see Section~\ref{sec:Cusp.Dirac}). Then, as shown in~\cite{Mo:MA08}, the Dirac Hamiltonian $(\shD_g)^2$ acting on the sections of the spinor bundle $\shS_X$ is essentially selfadjoint and has pure discrete spectrum. We then have an $n$-summable spectral triple, 
\begin{equation*}
 \left(\dot{C}^\infty(X), L^2_g(X,\shS_X), \shD_g\right). 
\end{equation*}
By using the Weyl laws of Moroianu~\cite{Mo:MA08} for fully elliptic cusp operators (which are derived from the analysis of their zeta functions), we show that we can apply to this spectral triple the semiclassical Weyl laws for spectral triples of~\cite{Po:Part1}. This gives the following semiclassical Weyl law (see Theorem~\ref{thm:cusp.SC-Weyl}).

\begin{theoremalph}\label{thm:Intro.SC-Weyl-Dirac}
Given any $q>0$ and any real-valued potential $V\in \dot{C}(X)$, for every energy level $ \lambda\in \R$, we have
\begin{gather*}
 \lim_{h\rightarrow 0^+} h^n N\left(h^{2q}\big(\shD_g^2\big)^{q}+V;\lambda\right) = \slashed{c}(n) \int_X (V(y)-\lambda)_{-}^{\frac{n}{2q}}d\nu_g(y), 
\end{gather*}
where $\slashed{c}(n):=2^{\left[\frac{n}{2}\right]}c(n)$ and $c(n)$ is the constant in~(\ref{eq:Intro.SC-Weyl-closed}). 
 \end{theoremalph}

\subsection{Magnetic Laplacians}
Let $\Delta_A=d_A^*d_A$ be the magnetic Laplacian, where $d_A=d+iA\wedge \cdot$ and $A$ is a cusp real $1$-form. Assume that $A$ is trapping in the sense of~\cite{GM:AHP08} (see Section~\ref{sec:Cusp.Magnetic}). Then, as shown in~\cite{GM:AHP08}, the operator $\Delta_A$ is essentially selfadjoint and has pure discrete spectrum. We then have an $n$-summable spectral triple, 
\begin{equation*}
 \left(\dot{C}^\infty(X), L^2_g(X), \sqrt{\Delta_A}\right). 
\end{equation*}
The results of Gol\'enia--Moroianu~\cite{GM:AHP08} allow us to apply the semiclassical Weyl laws of~\cite{Po:Part1} to this spectral triple as well. We then get the following result (see Theorem~\ref{thm:Cusp-SCWeyl-Magnetic}).

\begin{theoremalph}\label{thm:Intro.SC-Weyl-Magnetic}
Given any $q>0$ and any real-valued potential $V\in \dot{C}(X)$, for every energy level $ \lambda\in \R$, we have
\begin{gather*}
 \lim_{h\rightarrow 0^+} h^n N\left(h^{2q}\Delta_A^{q}+V;\lambda\right) = c(n) \int_X (V(y)-\lambda)_{-}^{\frac{n}{2q}}d\nu_g(y). 
\end{gather*}
 \end{theoremalph}

\subsection{Hodge Laplacians}
Let $\Delta_g^{(k)}$ be the Hodge Laplacian of $(X,g)$ acting on $k$-forms. Assume that the Betti numbers $b_{k-1}(M)$ and $b_k(M)$ both vanish (see Section~\ref{sec:Cusp.Hodge}). Then, as shown in~\cite{GM:TAMS12}, the operator $\Delta_g^{(k)}$ is essentially selfadjoint and has pure discrete spectrum. We then have an $n$-summable spectral triple, 
\begin{equation*}
 \left(\dot{C}^\infty(X), L^2_g(X,\Lambda^kT^*X), \sqrt{\Delta_g^{(k)}}\right). 
\end{equation*}
Here again, combining the results of Gol\'enia--Moroianu~\cite{GM:TAMS12} and Moroianu~\cite{Mo:MA08} with those of~\cite{Po:Part1}, we arrive at the following semiclassical Weyl law (see Theorem~\ref{thm:Cusp-SCWeyl-Hodge}).

\begin{theoremalph}\label{thm:Intro.SC-Weyl-Hodge}
Given any $q>0$ and any real-valued potential $V\in \dot{C}(X)$, for every energy level $ \lambda\in \R$, we have
\begin{gather*}
 \lim_{h\rightarrow 0^+} h^n N\Big(h^{2q}\big(\Delta_g^{(k)}\big)^{q}+V;\lambda\Big) = c(n,k) \int_X (V(y)-\lambda)_{-}^{\frac{n}{2q}}d\nu_g(y), \qquad c(n,k):=\binom{n}{k}c(n). 
\end{gather*}
 \end{theoremalph}

As mentioned above, the Riemannian manifold $(X,g)$ is not complete for $r>1$. Thus, in this case Theorems~\ref{thm:Intro.SC-Weyl-Dirac}--\ref{thm:Intro.SC-Weyl-Hodge} provide us with semiclassical Weyl laws on Riemannian manifolds that are not complete. For $r\leq 1$ they provide semiclassical Weyl laws on complete manifolds that do not have bounded geometry and with non-confining potentials (see Section~\ref{subsec:Intro.comparison}).

\subsection{Comparison with the results of Dai--Yan and Braverman}\label{subsec:Intro.comparison}
Dai--Yan~\cite{DY:arXiv25} established semiclassical Weyl laws for Schr\"odinger operators $h^2\Delta_g+V$ acting on functions on complete Riemannian manifolds with bounded geometry (i.e., the injectivity radius is bounded below and the curvature tensor and its first covariant derivative are bounded). They require the potential $V\in L^\infty_{\textup{loc}}$ to go to infinity at infinity and to satisfy a doubling condition and a regularity condition. They also obtain the classical Weyl law for $\Delta_g+V$. Subsequently, by using a localization argument, Braverman~\cite{Br:PAMS25} established semiclassical Weyl laws for the same type of operators on any \emph{complete} Riemannian manifold, where $V$ is a non-negative smooth function such that $V(x)\rightarrow \infty$ as $x\rightarrow \infty$. 

Our results lie outside the setting of~\cite{DY:arXiv25, Br:PAMS25} in the following respects. First, the manifold $(X,g)$ is not complete for $r>1$ (see~\cite{GM:TAMS12, Mo:PC}), and so in this case it falls outside the scope of~\cite{DY:arXiv25, Br:PAMS25}. 

Second, for $r\leq 1$ the manifold $(X,g)$ is complete, but it does not have bounded geometry, and so it falls outside the scope of~\cite{DY:arXiv25}. Indeed, as $(X,g)$ is complete, non-compact, and has finite volume, it contains infinitely many disjoint geodesic balls of radius~$1$, and the volumes of these balls must go to $0$. However, bounded geometry implies a uniform positive lower bound on the volumes of geodesic balls of radius~$1$ (see also~\cite[p.~140]{GM:AHP08}, where it is observed that the injectivity radius goes to $0$ at infinity). 

Third, in~\cite{DY:arXiv25, Br:PAMS25} it is the confining property of the potential that ensures that the Schr\"odinger operators $h^2\Delta_g+V$ have discrete spectra. In contrast, in our setting the operators $\shD_g^2$, $\Delta_A$, and $\Delta_g^{(k)}$ have discrete spectra by themselves, and the potentials are bounded. In particular, for $r\leq 1$ our results are not covered by~\cite{Br:PAMS25}, even though the manifold is complete. 

Finally, our results hold for fractional powers of the operators and for operators acting on bundles, whereas~\cite{DY:arXiv25, Br:PAMS25} consider Schr\"odinger operators $h^2\Delta_g+V$ acting on functions. 

\subsection{Integration formulas}
In the framework of noncommutative geometry, the role of the integral is played by the noncommutative integral $\bint$ of Connes~\cite{Co:NCG} (see Section~\ref{sec:Part1.Integration}). If $(M^n,g)$ is a closed Riemannian manifold, then Connes' integration formula asserts that, for every $f\in C(M)$, we have
\begin{equation*}
 \bint f\Delta_g^{-\frac{n}{2}}= c(n) \int_M f(x) d\nu_g(x).
\end{equation*}
Thus, the noncommutative integral recaptures the Riemannian measure. The main results of~\cite{Po:Part1} include extensions of this formula to spectral triples. In the setting of this paper they lead to the following integration formula, which is a direct extension to open manifolds of the original integration formula of Connes.

\begin{theoremalph}\label{thm:Intro.Integration-Dirac}
 For every $f\in  \dot{C}(X)$, the operator $f|\shD_g|^{-n}$ is spectrally measurable, with
  \begin{gather*}
 \bint f|\shD_g|^{-n}=\slashed{c}(n)\int_X f(y) d\nu_g(y).
\end{gather*}
There are similar results for the absolute values $|f|\shD_g|^{-n}|$, and for the operators $|\shD_g|^{-n/2}f|\shD_g|^{-n/2}$ and their absolute values.
\end{theoremalph}

Here again, there are similar integration formulas for magnetic Laplacians and Hodge Laplacians (see Theorem~\ref{thm:Cusp.Integration-Laplacians}).

\subsection{Conformally cusp metrics of infinite volume}
The condition $r>1/n$ ensures that $(X,g)$ has finite volume. For $0<r\leq 1/n$ the manifold $(X,g)$ is complete and has infinite volume, and the results of~\cite{Mo:MA08, GM:AHP08, GM:TAMS12} still provide us with Weyl laws for the above operators. For $r<1/n$ these Weyl laws take the form $\lambda_j \sim (j/\gamma)^{2r}$, where the constant $\gamma$ only depends on the geometry at infinity. In this case the results of~\cite{Po:Part1} apply as well, but the semiclassical Weyl laws and the integration formulas only involve the values at infinity. For instance, in the case of the Dirac operator we obtain the following result (see Theorem~\ref{thm:Cusp.infinite}).

\begin{theoremalph}\label{thm:Intro.SC-Weyl-infinite}
 Suppose that $0<r<1/n$, and assume that $\overline{X}$ is spin, $g_0$ is an exact cusp metric, and there are no non-zero harmonic spinors on $M$. Given any $q>0$ and any real-valued potential $V\in \dot{C}(X)$, for every energy level $ \lambda\in \R$, we have
\begin{equation*}
 \lim_{h\rightarrow 0^+} h^{\frac{1}{r}} N\left(h^{2q}\big(\shD_g^2\big)^{q}+V;\lambda\right) = \gamma \left(V_0-\lambda\right)_{-}^{\frac{1}{2qr}},
\end{equation*}
where $V_0$ is the boundary value of $V$ and $\gamma$ is a positive constant which only depends on $r$ and the boundary $(M,h_0)$ (see Eq.~(\ref{eq:Cusp.infinite-gamma})). 
\end{theoremalph}

There are similar results for magnetic Laplacians and Hodge Laplacians, as well as integration formulas (see Theorem~\ref{thm:Cusp.infinite}). For $r=1/n$ the Weyl laws are nonclassical, and so the results of~\cite{Po:Part1} do not apply (see Section~\ref{subsec:Cusp.critical}). 

\subsection{Organization of the paper}
The rest of the paper is organized as follows. In Section~\ref{sec:Part1.SC}, we recall the main results of~\cite{Po:Part1} on semiclassical Weyl laws for spectral triples. In Section~\ref{sec:Cusp.metrics}, we recall the main definitions regarding conformally cusp metrics. In Section~\ref{sec:Cusp.Dirac}, we construct the spectral triple associated with the Dirac operator and establish the semiclassical Weyl laws for this operator. In Section~\ref{sec:Cusp.Magnetic}, we establish the corresponding results for magnetic Laplacians. In Section~\ref{sec:Cusp.Hodge}, we obtain semiclassical Weyl laws for Hodge Laplacians. In Section~\ref{sec:Part1.Integration}, we review Connes' integration and the integration formulas of~\cite{Po:Part1}. In Section~\ref{sec:Cusp.Integration}, we derive the integration formulas for Dirac operators, magnetic Laplacians, and Hodge Laplacians. In Section~\ref{sec:Cusp.infinite}, we deal with conformally cusp metrics of infinite volume.

\subsection*{Acknowledgements}
The author wishes to thank Sergiu Moroianu for various discussions related to the subject matter of this article. The research for this article was partially funded by NSFC grant No.~11971328 (China).

\section{Semiclassical Weyl Laws for Spectral Triples}\label{sec:Part1.SC}
In this section, we recall the main results of~\cite{Po:Part1} on semiclassical Weyl laws for spectral triples.

\subsection{Spectral triples}
In the framework of noncommutative geometry, noncommutative manifolds are represented by \emph{spectral triples} $(\sA,\sH,D)$, where $\sA$ is a (unital) $*$-algebra represented by bounded operators on the Hilbert space $\sH$ and $D$ is an unbounded selfadjoint operator on $\sH$ with compact resolvent such that the commutators $[D,a]$, $a\in \sA$, are bounded. The notion of dimension is replaced by $p$-summability, in the sense that $D^{-1}$ lies in the weak Schatten class $\sL_{p,\infty}$. That is, 
\begin{equation*}
 \mu_j(D^{-1})=\op{O}\big(j^{-\frac{1}{p}}\big) \qquad \text{as}\ j\rightarrow \infty, 
\end{equation*}
where $\mu_0(D^{-1})\geq \mu_1(D^{-1}) \geq \cdots $ are the singular values of $D^{-1}$, i.e., the eigenvalues of its absolute value $|D^{-1}|$.

The main prototype of a spectral triple is the Dirac spectral triple $(C^\infty(M), L^2(M,\shS),\shD_g)$ associated with any closed Riemannian spin manifold $(M^n,g)$. A closely related example is the square-root Laplacian spectral triple $(C^\infty(M),L^2(M),\sqrt{\Delta_g})$. In both cases we have $n$-summable spectral triples with $n=\dim M$.

\subsection{Semiclassical Weyl laws}
The spectral asymptotics for eigenvalues of Birman--Schwinger operators form the backbone of the operator-theoretic approach to semiclassical analysis of Birman--Solomyak~\cite{BS:FAA70, BS:AMST89} (see also Simon~\cite{Si:TAMS76}). In the case of a closed Riemannian manifold $(M^n,g)$, for any real-valued potential $V(x)$ which is continuous (or even in some $L^p$-space), the Birman--Schwinger principle relates Weyl laws for the Birman--Schwinger operators $K_V^{(q)}(\lambda):=\Delta_g^{-q/2}(V-\lambda)\Delta_g^{-q/2}$, $q>0$, $\lambda \in \R$, to semiclassical Weyl laws for fractional Schr\"odinger operators $H_V^{(q)}(h):=h^{2q}\Delta_g^{q}+V$. Namely, if $N(H_V^{(q)}(h);\lambda)$, $\lambda\in \R$, is the counting function of $H_V^{(q)}(h)$, then
\begin{align*}
 \lim_{h\rightarrow 0^+} h^n N\big(H_V^{(q)}(h);\lambda\big) = \lim_{j\rightarrow \infty} j\lambda_j^{-}\big(K_V^{(q)}(\lambda)\big)^{\frac{n}{2q}}
 = c(n) \int_M \left(V(x)-\lambda\right)_{-}^{\frac{n}{2q}} d\nu_g(x). 
\end{align*}
The first equality is the Birman--Schwinger principle. The second equality follows from the Weyl laws of Birman--Solomyak~\cite{BS:VLU77, BS:VLU79, BS:SMJ79} for negative order \psidos. For $q=1$ this recovers the semiclassical Weyl law of Colin de Verdi\`ere~\cite{CdV:DM79} without remainder term.

Suppose that $(\sA, \sH, D)$ is a $p$-summable spectral triple. Let $\overline{\sA}$ be the closure of $\sA$ in $\sL(\sH)$ (this is a $C^*$-algebra). Denote by 
$\sA_{++}$ the cone of positive invertible elements of $\sA$. For every $a\in \sA_{++}$, the operator $aD^2a$ has a non-negative discrete spectrum, and so we can arrange its eigenvalues as a non-decreasing sequence, 
\begin{equation*}
 \lambda_0\left(aD^2a\right) \leq \lambda_1\left(aD^2a\right) \leq \lambda_2\left(aD^2a\right) \leq \cdots,
\end{equation*}
where each eigenvalue is repeated according to multiplicity. We assume that, for all $a\in \sA_{++}$, we have
\begin{equation}\label{eq:Intro.conf-Weyl}
 \lim_{j\rightarrow \infty}j^{-\frac{2}{p}}\lambda_j\big(aD^2a\big)= \tau\big[a^{-p}\big]^{-\frac{2}{p}}, 
\end{equation}
where $\tau:\overline{\sA}\rightarrow \C$ is a given positive linear map (which does not depend on $a$). This spectral condition is called \emph{Condition~(W)} in~\cite{Po:Part1}. This amounts to requiring to have Weyl laws for conformal deformations of the square $D^2$, where conformal deformation is meant in the sense of Connes--Moscovici~\cite{CM:Twisted}.

\begin{theorem}[{\cite[Theorem~B]{Po:Part1}}]\label{thm:Part1.SC-Weyl}
 Assume that Condition~\textup{(W)} is satisfied. Given any $q>0$ and  $V^*=V\in \overline{\sA}$, for every energy level $ \lambda\in \R$, we have
\begin{equation}\label{eq:Part1.SC-Weyl}
 \lim_{h\rightarrow 0^+} h^{p} N\big(h^{2q}(D^2)^{q}+V;\lambda\big)= \tau\Big[\left(V-\lambda\right)_{-}^{\frac{p}{2q}}\Big].
\end{equation}
\end{theorem}

Unlike in~\cite{MSZ:LMP22}, the above theorem holds for any summability degree $p>0$ and without extra regularity assumptions on the spectral triple. By the Birman--Schwinger principle it is a consequence of Weyl laws for the Birman--Schwinger operators $|D|^{-q}(V-\lambda)|D|^{-q}$, which are established in~\cite{Po:Part1} as well (see Section~\ref{sec:Part1.Integration}).

Applying Theorem~\ref{thm:Part1.SC-Weyl} requires two steps:
\begin{enumerate}
 \item Construct a spectral triple $(\sA, \sH,D)$. In particular, it is crucial to establish boundedness of the commutators $[D,a]$, $a\in \sA$.

 \item Verify that Condition~(W) is satisfied.
\end{enumerate}

\section{Conformally Cusp Metrics}\label{sec:Cusp.metrics}
Throughout this paper, we let $X^n$ be an open manifold, which is the interior of a compact manifold $\overline{X}$ with closed boundary $M^{n-1}=\partial \overline{X}$. Let $x$ be a boundary defining function, i.e., $x$ is a non-negative smooth function on $\overline{X}$ such that $x^{-1}(0)=M$ and $dx\neq 0$ on $M$. We assume that $X$ is endowed with a Riemannian metric of the form, 
\begin{equation*}
 g=x^{2r}g_0, \qquad r>0,
\end{equation*}
where $g_0$ is a cusp metric on $\overline{X}$, i.e., a (complete) Riemannian metric which in local coordinates near $M$ takes the form, 
\begin{equation*}
 g_0= a_{00}(x,y)x^{-4}dx^2+ \sum_{1\leq j<n}a_{0j}(x,y)x^{-2}dxdy^{j} +  \sum_{1\leq i, j<n}h_{ij}(x,y)dy^i dy^j. 
\end{equation*}
 If in addition $a_{00|M}=1$ and $a_{0j|M}=0$ for $j=1, \ldots, n-1$, then (following Melrose~\cite{Me:Book}) we say that $g_0$ is an \emph{exact} cusp metric. Note that the cusp metric $g_0$ induces on $M$ the metric, 
\begin{equation}\label{eq:Cusp.boundary-metric}
 h_0:=  \sum_{1\leq i, j<n}h_{ij}(0,y)dy^i dy^j. 
\end{equation}
The condition $r>1/n$ is exactly the condition that ensures that $(X,g)$ has finite volume. In addition, $(X,g)$ is a complete Riemannian manifold if and only if $r\leq 1$ (see~\cite{GM:TAMS12, Mo:PC}). Except in Section~\ref{sec:Cusp.infinite}, where we deal with the case $0<r\leq 1/n$, we assume that $r>1/n$ throughout the paper. 

Examples of manifolds with conformally cusp metrics  with $r> 1/n$ include hyperbolic metrics of finite volume (in which case $r=1$; see~\cite{Mo:MA08}).  Another class of examples is provided by metric horns in the sense of Cheeger~\cite{Ch:PSPM80} (see also~\cite{LP:CCPDE98}) for which $r>1$. In particular, these metrics are incomplete.

\section{Semiclassical Weyl Laws for Dirac Operators}\label{sec:Cusp.Dirac} Assume further that $\overline{X}$ is spin. We denote by  $\shD_g$ the Dirac operator of $X$ acting on the sections of the spinor bundle $\shS_X$. In addition, we denote by $\shD_{h_0}$ the Dirac operator of $M$ with respect to the metric~(\ref{eq:Cusp.boundary-metric}).  We also make the following assumptions: 
\begin{itemize}
 \item $g_0$ is an exact cusp metric. 
 
 \item $\ker \shD_{h_0}=\{0\}$, i.e., there are no (non-zero) harmonic spinors on $M$. 
\end{itemize}
The second condition holds for instance if $(M,h_0)$ has a non-negative scalar curvature which is non-zero at a point.  It is also known to hold for generic metrics if $M$ either has dimension $\leq 4$ or its fundamental group $\pi_1(M)$ vanishes or is cyclic of odd order (\emph{cf}.\ the discussion in~\cite[p.~3]{Mo:MA08} and the references therein). 

We have a natural unitary isomorphism $U:L_{g_0}^2(X,\shS_X)\rightarrow L_g^2(X,\shS_X)$ given by
\begin{equation}\label{eq:Cusp.isometry}
 U\phi = x^{-\frac{nr}{2}}\phi, \qquad \phi \in L^2_{g_0}(X,\shS_X). 
\end{equation}
As observed in~\cite{Mo:MA08}, under these assumptions the conjugate operator $\hat{D}_g:=U^*\shD_gU$ is a fully elliptic cusp differential operator of type $(1,r)$. Full ellipticity means that the principal symbol of $x^r\hat{D}_g$ is invertible and the normal operator
\begin{equation*}
 \sN(x^r\hat{D}_g)(\xi):=(e^{i\xi/x}x^r\hat{D}_ge^{-i\xi/x})_{|M}, \qquad \xi\in \R,
\end{equation*}
is invertible for all $\xi\in \R$. This ensures the existence of a parametrix for $\hat{D}_g^2$ in the cusp pseudodifferential calculus of Melrose (see, e.g.,~\cite[Appendix~A]{MN:arXiv96}). Using this \psido\ calculus the following results are established.  
\begin{enumerate}
\item[(i)]  $\shD_g$ with domain $C_c^\infty(X,\shS_X)$ is essentially selfadjoint on $L_g^2(X,\shS_X)$ and has compact resolvent. Thus, its spectrum consists of  isolated real eigenvalues with finite multiplicity (see also~\cite{Ba:JDG00} for the case $r=1$).  

\item[(ii)] For every $z\in \C$, the power $|\shD_g|^{z}$ is such that $U^*|\shD_g|^{z}U=|\hat{D}_g|^{z}$ is a cusp \psido\ of type $(z,rz)$. In particular, it is trace-class for 
$\Re z<-n$. 

\item[(iii)] The zeta function $\zeta_{\shD_g}(z):=\Tr[|\shD_g|^{-z}]$, $\Re z>n$, has a meromorphic continuation to $\C$ with at worst poles on $(n-\N_0)\cup (r^{-1}-r^{-1}\N_0)$. 

\item[(iv)] The first pole occurs at $z=n$. This pole is simple and the residue is given by
\begin{equation*}
 \Res_{z=n}\zeta_{\shD_g}(z)= n\slashed{c}(n)\Vol_g(X), \qquad \text{where}\ \slashed{c}(n):= 2^{\left[\frac{n}{2}\right]}(2\pi)^{-n}|\bB^{n}|.
\end{equation*}

\item[(v)] By Ikehara's Tauberian theorem this implies that, as $j\rightarrow \infty$, we have the Weyl law, 
\begin{equation*}
 \lambda_j\big(\shD_g^2\big) \sim j^{\frac{2}{n}} \big(\slashed{c}(n)\Vol_g(X)\big)^{-\frac{2}{n}}. 
 \end{equation*}
\end{enumerate}

Bearing this in mind, to get a spectral triple we use the $*$-algebra, 
\begin{equation*}
 \dot{C}^\infty(X)=\C1+C_c^\infty(X). 
\end{equation*}
That is, $\dot{C}^\infty(X)$ is the algebra of $C^\infty$-functions on $X$ that are constant near the boundary $M=\partial \overline{X}$.  It is represented in $L_g^2(X,\shS_X)$ by bounded multiplication operators. Its norm closure is the $C^*$-algebra,
\begin{equation*}
 \dot{C}(X):=\C1+C_0(X),
\end{equation*}
 where $C_0(X)$ is the ideal of continuous functions on $X$ that converge to $0$ near the boundary $M$. Thus, $\dot{C}(X)$ consists of continuous functions on $\overline{X}$ that are constant on the boundary.  If $f$ is such a function, we denote by $f_0\in \C$ its boundary value. 

If $f\in \dot{C}^\infty(X)$, then $[\shD_g,f]=c(df)$ is a compactly supported section of $\End \shS_X$, and hence it is bounded on $L_g^2(X,\shS_X)$. We thus get a spectral triple, 
\begin{equation}\label{eq:Cusp.Dirac-ST}
 \left(\dot{C}^\infty(X), L^2_g(X,\shS_X), \shD_g\right). 
\end{equation}
This spectral triple is $n$-summable.

Let $f\in \dot{C}^\infty(X)$ be~$>0$ on $X$ and have positive boundary value. The multiplication by $f$ commutes with the unitary isomorphism $U$ given by~(\ref{eq:Cusp.isometry}). Thus,
\begin{equation*}
 U^*(f\shD_g^2f)U=f(U^*\shD_gU)^2f =f\hat{D}_g^2f.  
\end{equation*}
 The operator $f\hat{D}_g^2f$ is a fully elliptic cusp differential operator of type $(2,2r)$, since its principal symbol $f^2|\xi|_{g_0}^2$ and its normal operator $f_0^2\sN(x^{r}\hat{D}_g)^2$ are both positive. Therefore, the results of~\cite{Mo:MA08} apply, and so as  $j\rightarrow \infty$ we have the Weyl law, 
\begin{equation*}
 \lambda_j\big(f\shD_g^2f\big) \sim  
j^{\frac{2}{n}} \left(\slashed{c}(n)\int_X f(y)^{-n} d\nu_g(y)\right)^{-\frac{2}{n}}.
\end{equation*}

All this shows that $(\dot{C}^\infty(X), L^2_g(X,\shS_X), \shD_g)$ is a spectral triple satisfying Condition~\textup{(W)} with
\begin{equation}\label{eq:Cusp.tau-Dirac}
 \tau(f)=\slashed{c}(n)\int_X f(y) d\nu_g(y), \qquad f\in \dot{C}(X). 
\end{equation}
Therefore, Theorem~\ref{thm:Part1.SC-Weyl} yields the following semiclassical Weyl law. 

\begin{theorem}\label{thm:cusp.SC-Weyl}
 Let $q>0$, and assume that $g_0$ is an exact cusp metric and there are no non-zero harmonic spinors on $M$. 
 Given any real-valued potential $V\in \dot{C}(X)$, for every energy level $\lambda \in \R$, we have
\begin{gather*}
 \lim_{h\rightarrow 0^+} h^n N\left(h^{2q}\big(\shD_g^2\big)^{q}+V;\lambda\right) = \slashed{c}(n) \int_X (V(y)-\lambda)_{-}^{\frac{n}{2q}}d\nu_g(y).
\end{gather*}
\end{theorem}

\begin{remark}
 As pointed out to the author by Sergiu Moroianu~\cite{Mo:PC}, in the cases $r=1,2$ we can replace $\dot{C}^\infty(X)$ by a larger algebra to get a spectral triple in~(\ref{eq:Cusp.Dirac-ST}). For $r=1$ (i.e., geometric cusp case) we may use the algebra of functions that are smooth up to the boundary and constant on $M$. For $r=2$ we need to further require the normal derivative (i.e., $\partial f/\partial x$) is constant on $M$ as well. In either case we get the same norm closure $\dot{C}(X)$ as with using $\dot{C}^\infty(X)$, and so this does not affect our main results. 
\end{remark}

\section{Semiclassical Weyl Laws for Magnetic Laplacians}\label{sec:Cusp.Magnetic} We drop  the assumption on $X$ being spin. The ordinary Laplacian of $(X,g)$ need not have discrete spectrum (see, e.g., \cite{FH:DMJ89}). However, as observed by Gol\'enia--Moroianu~\cite{GM:AHP08} we get a fully elliptic cusp differential operator with discrete spectrum by coupling the Laplacian with a trapping magnetic field. Such operators appear in the study of the Aharonov--Bohm effect. 

More precisely, let $B$ be  a magnetic field of the form $B=dA$, where the vector potential $A$ is a  cusp real 1-form, i.e., a smooth real 1-form on $X$ which near $M$ is of the form, 
\begin{equation}\label{eq:cusp.cusp-1-form}
 A=a_{0}(x,y) x^{-2}dx +\omega, \qquad \iota_{\partial_x}\omega=0, 
\end{equation}
where $a_0(x,y)$ and $\omega$ are smooth up to the boundary. Set $A_0=\omega_{|M}$. Following~\cite{GM:AHP08} we shall say that $A$ is \emph{trapping} on a given connected component $M_0\subseteq M$ if one of the following two conditions is satisfied:
\begin{itemize}
 \item $a_0(0,y)$ is not constant on $M_0$. 
 
 \item Either $A_0$ is not closed on $M_0$, or it is closed and its cohomology class in $H^1(M_0,\R)$ does not lie in $2\pi H^1(M_0;\Z)$, i.e., $(2\pi)^{-1}A_0$ defines a non-integral cohomology class. 
\end{itemize}
We say that $A$ is trapping on $M$ if it is trapping on every connected component. 

The \emph{magnetic Laplacian} is the differential operator on $X$ defined by
\begin{equation*}
 \Delta_{A}:=d_A^*d_A, \qquad d_A=d+iA\wedge \cdot,
\end{equation*}
where $d:C_c^\infty(X)\rightarrow C_c^\infty(X,T^*X)$ is the de Rham differential. 

As shown in~\cite{GM:AHP08}, if $A$ is trapping, then $\Delta_A$ is a fully elliptic cusp differential operator of type $(2,2r)$. Thus, in the same way as with the square of the Dirac operator,  the operator $\Delta_A$ with domain $C_c^\infty(X)$ is essentially selfadjoint on $L_g^2(X)$ and has compact resolvent. Moreover, we have the Weyl law, 
\begin{equation*}
 \lambda_j\big(\Delta_A\big) \sim  j^{\frac{2}{n}} \big(c(n)\Vol_g(X)\big)^{-\frac{2}{n}}, \qquad c(n):= (2\pi)^{-n}|\bB^{n}|. 
 \end{equation*}

If $f\in \dot{C}^\infty(X)$, then $[\sqrt{\Delta_A},f]=[\sqrt{\Delta_A},f-f_0]$. Here $f-f_0\in C^\infty_c(X)$, and so $\sqrt{\Delta_A}(f-f_0)$ and 
 $(f-f_0)\sqrt{\Delta_A}$ are cusp \psidos\ of type $(1,-\infty)$ with the same principal symbol and vanishing normal operators. It follows that $[\sqrt{\Delta_A},f]$
 is a cusp \psido\ of type $(0,-\infty)$, and hence is bounded on $L_g^2(X)$. We thus have a spectral triple, 
\begin{equation}\label{eq:Cusp.Magnetic-ST}
 \left(\dot{C}^\infty(X), L^2_g(X), \sqrt{\Delta_A}\right). 
\end{equation}
As with the Dirac spectral triple~(\ref{eq:Cusp.Dirac-ST}), this spectral triple is $n$-summable.

If $f\in  \dot{C}^\infty(X)$ is~$>0$ on $X$ and has a positive boundary value, then $f\Delta_Af$ is again a fully elliptic cusp differential operator of type $(2,2r)$. Indeed, the principal symbol of $x^{2r}f\Delta_Af$ is $f^2|\xi|^2_{g_0}$ and its normal operator is $\sN(x^{2r}f\Delta_Af)(\xi)=f_0^2\sN(x^{2r}\Delta_A)(\xi)$. Therefore, in the same way as with $\shD^2_g$, as $j\rightarrow \infty$ we have the Weyl law, 
\begin{equation}\label{eq:Cusp.Weyl-Magnetic}
 \lambda_j\big(f\Delta_A f\big) \sim j^{\frac{2}{n}} \left(c(n)\int_X f(y)^{-n} d\nu_g(y)\right)^{-\frac{2}{n}}. 
 \end{equation}

It follows from all this that $(\dot{C}^\infty(X), L^2_g(X), \sqrt{\Delta_A})$ is a spectral triple satisfying Condition~\textup{(W)} with
\begin{equation}\label{eq:Cusp.tau-Magnetic}
 \tau(f)=c(n)\int_X f(y) d\nu_g(y), \qquad f\in \dot{C}(X). 
\end{equation}
Therefore, in the same way as with the Dirac operator, Theorem~\ref{thm:Part1.SC-Weyl} yields the following semiclassical Weyl law. 
 
\begin{theorem}\label{thm:Cusp-SCWeyl-Magnetic}
 Let $q>0$, and assume that the $1$-form $A$ is trapping. For any real-valued potential $V\in \dot{C}(X)$ and every energy level $\lambda \in \R$, we have 
 \begin{gather}\label{eq:Cusp.SCWeyl-Magnetic1}
 \lim_{h\rightarrow 0^+} h^n N\left(h^{2q}\Delta_A^q+V;\lambda\right) = c(n) \int_X (V(y)-\lambda)_{-}^{\frac{n}{2q}}d\nu_g(y).
\end{gather}
\end{theorem}
 
 \begin{remark}
 We refer to~\cite{MOR:AHP05, RS:AMS99} for Cwikel--Lieb--Rozenblum inequalities and semiclassical Weyl laws for magnetic Schr\"odinger operators on $\R^2$ and $\R^3$. 
\end{remark}
 
\begin{remark}
 We also refer to~\cite{LS:JAM97} for a semiclassical Weyl law for Schr\"odinger operators $h^2\Delta_{\text{hyp}}+V$ on hyperbolic space $\mathbb{H}^{n}$, $n\geq 3$. 
\end{remark}

\section{Semiclassical Weyl Laws for Hodge Laplacians}\label{sec:Cusp.Hodge}
In this section, we explain how the results of Gol\'enia--Moroianu~\cite{GM:TAMS12} allow us to get semiclassical Weyl laws for Hodge Laplacians on $(X,g)$. As in Section~\ref{sec:Cusp.Magnetic}, we do not assume $X$ to be spin, nor do we assume $g_0$ to be an exact cusp metric. 

Given $k\in \{0,1,\ldots,n\}$, let $\Delta_g^{(k)}=d^*d+dd^*$ be the Hodge Laplacian of $(X,g)$ acting on $k$-forms with domain $C^\infty_c(X,\Lambda^kT^*X)$. Denote by $b_j(M)$, $j=0,1,\ldots, n-1$, the Betti numbers of $M$. As shown in~\cite[Proposition~3.1]{GM:TAMS12}, the operator $\Delta_g^{(k)}$ belongs to $x^{-2r}\op{Diff}^2_c(X,\Lambda^k({}^cT X))$, it is cusp-elliptic, and $x^{2r}\Delta_g^{(k)}$ is fully elliptic if and only if $b_{k-1}(M)=b_k(M)=0$. 

From now on we assume that $b_{k-1}(M)=b_k(M)=0$. Note that this condition does not hold for $k=0,1$. It does not hold for $k=n-1,n$ either if $M$ has at least one orientable connected component (see~\cite{GM:TAMS12}). Under this assumption, the following results are established in~\cite[Corollary~3.2 and Theorem~1.2]{GM:TAMS12}. 
\begin{enumerate}
\item[(i)] The operator $\Delta_g^{(k)}$ is essentially selfadjoint on $L_g^2(X,\Lambda^kT^*X)$ and has pure discrete spectrum. Moreover, the domain of its selfadjoint extension is the weighted cusp Sobolev space $x^{2r}H^2_c(X,\Lambda^kT^*X)$ defined in~\cite{GM:TAMS12}. 

\item[(ii)] The complex powers of $\Delta_g^{(k)}$ are cusp \psidos, and the zeta function $\Tr[(\Delta_g^{(k)})^{-z/2}]$, $\Re z>n$, has a meromorphic continuation to $\C$. As $r>1/n$, its first pole is simple and located at $z=n$, with residue
\begin{equation*}
 (2\pi)^{-n}\binom{n}{k}\left|\bS^{n-1}\right|\Vol_g(X)=n c(n,k)\Vol_g(X), 
\end{equation*}
where we have set $c(n,k):=\binom{n}{k}c(n)=(2\pi)^{-n}\binom{n}{k}| \bB^n|$. 

\item[(iii)] As $j\rightarrow \infty$, we have the Weyl law, 
\begin{equation*}
 \lambda_j\big(\Delta_g^{(k)}\big) \sim  j^{\frac{2}{n}} \big(c(n,k)\Vol_g(X)\big)^{-\frac{2}{n}}. 
 \end{equation*}
\end{enumerate}

If $f\in \dot{C}^\infty(X)$, then, as $\Delta_g^{(k)}$ has a scalar principal symbol, in the same way as with the magnetic Laplacian, we see that $[\sqrt{\Delta_g^{(k)}},f]=[\sqrt{\Delta_g^{(k)}},f-f_0]$ is a cusp \psido\ of type $(0,-\infty)$, and hence is bounded on $L_g^2(X,\Lambda^kT^*X)$. We thus have a spectral triple, 
\begin{equation}\label{eq:Cusp.Hodge-ST}
 \left(\dot{C}^\infty(X), L^2_g(X,\Lambda^kT^*X), \sqrt{\Delta_g^{(k)}}\right). 
\end{equation}
This spectral triple is $n$-summable. 

Let $f\in  \dot{C}^\infty(X)$ be~$>0$ on $X$ and have a positive boundary value $f_0$. Then $f\Delta_g^{(k)}f$ is a symmetric cusp-elliptic operator in $x^{-2r}\op{Diff}^2_c(X,\Lambda^k({}^cT X))$, and the principal symbol of $x^{2r}f\Delta_g^{(k)}f$ is $f^2|\xi|_{g_0}^2$. Moreover, the normal operator of $x^{2r}f\Delta_g^{(k)}f$ is $f_0^2\sN(x^{2r}\Delta_g^{(k)})$, and so $x^{2r}f\Delta_g^{(k)}f$ is fully elliptic. Therefore, as in the case of the magnetic Laplacian, the results of~\cite{GM:TAMS12, Mo:MA08} apply, and so as $j\rightarrow \infty$ we have the Weyl law, 
\begin{equation*}
 \lambda_j\big(f\Delta_g^{(k)} f\big) \sim j^{\frac{2}{n}} \left(c(n,k)\int_X f(y)^{-n} d\nu_g(y)\right)^{-\frac{2}{n}}. 
 \end{equation*}

It follows that the spectral triple~(\ref{eq:Cusp.Hodge-ST}) satisfies Condition~\textup{(W)} with
\begin{equation}\label{eq:Cusp.tau-Hodge}
 \tau(f)=c(n,k)\int_X f(y) d\nu_g(y), \qquad f\in \dot{C}(X). 
\end{equation}
Therefore, Theorem~\ref{thm:Part1.SC-Weyl} yields the following semiclassical Weyl law. 

\begin{theorem}\label{thm:Cusp-SCWeyl-Hodge}
 Let $q>0$ and $k\in \{0,1,\ldots,n\}$, and assume that $b_{k-1}(M)=b_k(M)=0$. For any real-valued potential $V\in \dot{C}(X)$ and every energy level $\lambda \in \R$, we have 
 \begin{gather*}
 \lim_{h\rightarrow 0^+} h^n N\Big(h^{2q}\big(\Delta_g^{(k)}\big)^q+V;\lambda\Big) = c(n,k) \int_X (V(y)-\lambda)_{-}^{\frac{n}{2q}}d\nu_g(y).
\end{gather*}
\end{theorem}

\begin{remark}
 As observed in~\cite{GM:TAMS12}, the assumption $b_{k-1}(M)=b_k(M)=0$ ensures that $\Delta_g^{(k)}$ is essentially selfadjoint and has pure discrete spectrum regardless of whether the metric $g$ is complete. In particular, for $r>1$ Theorem~\ref{thm:Cusp-SCWeyl-Hodge} provides us with semiclassical Weyl laws for Hodge Laplacians on Riemannian manifolds that are not complete. 
\end{remark}

\section{Connes' Integration Formulas for Spectral Triples}\label{sec:Part1.Integration}
In this section, we review Connes' integration and the main results of~\cite{Po:Part1} on integration formulas for spectral triples.

\subsection{Connes' integration}
In the framework of noncommutative geometry, the role of the integral is played by (normalized) positive traces  on the weak trace class $\sL_{1,\infty}$, and more especially Dixmier traces (see~\cite{Co:NCG}). An operator $A\in \sL_{1,\infty}$ is then called \emph{measurable} if it takes the same value on all Dixmier traces. Equivalently, $A$ is measurable if and only if
\begin{equation*}
 \bint A := \lim_{N\rightarrow \infty}  \frac{1}{\log N}\sum_{j<N} \lambda_j(A)\ \text{exists}, 
\end{equation*}
where $\{\lambda_j(A)\}$ is an eigenvalue sequence for $A$ (see~\cite{LSZ:Book, Po:JNCG23}). The limit $\bint A$ is precisely the NC integral of $A$. A stronger notion of measurability (called \emph{strong measurability}) further requires $A$ to take the same value on \emph{all normalized positive traces}.

For instance, if $(M^n,g)$ is a closed Riemannian manifold with (positive) Laplace operator $\Delta_g$, then Connes' integration formula asserts that, for every $f\in C(M)$, the operator $f\Delta_g^{-n/2}$ is strongly measurable, with 
\begin{equation*}
 \bint f\Delta_g^{-\frac{n}{2}}= c(n) \int_M f(x) d\nu_g(x), 
\end{equation*}
where $c(n)$ is the constant in~(\ref{eq:Intro.SC-Weyl-closed}) and $\nu_g(x)$ is the Riemannian measure. This shows that the NC integral recaptures the Riemannian measure.

This notion of strong measurability is closely related to Weyl laws for eigenvalues of operators in weak Schatten classes. Namely, if $T=T^*\in \sL_{1,\infty}$, we say that $T$ is \emph{spectrally measurable} if the limits $\lim_{j\rightarrow \infty} j\lambda_j^\pm(T)$ both exist, in which case we have 
\begin{equation*}
 \bint T = \lim_{j \rightarrow \infty} j\lambda_j^+(T)- \lim_{j \rightarrow \infty} j\lambda_j^{-}(T). 
\end{equation*}
Spectral measurability implies strong measurability (see, e.g., \cite{Po:JNCG23}), but the converse does not hold (see, e.g., \cite{PT:Part2}). 

\subsection{Integration formulas for spectral triples}
Let $(\sA, \sH, D)$ be a $p$-summable spectral triple satisfying Condition~(W) (see Section~\ref{sec:Part1.SC}). It is shown in~\cite{Po:Part1} that, given any $q>0$, for every $a\in \overline{\sA}$, we have
\begin{equation}\label{eq:Part1.Spectral-asymp1}
 \lim_{j\rightarrow \infty}j^{\frac{q}{p}} \mu_j\big(|D|^{-\frac{q}{2}}a|D|^{-\frac{q}{2}}\big) = \lim_{j\rightarrow \infty}j^{\frac{q}{p}} \mu_j\big(a|D|^{-q}\big) = \tau\left[|a|^{\frac{p}{q}}\right]^{\frac{q}{p}}.
\end{equation}
If in addition $a$ is selfadjoint, then we have the Weyl laws,
\begin{equation}\label{eq:Part1.Spectral-asymp2}
 \lim_{j\rightarrow \infty} j^{\frac{q}{p}} \lambda_j^\pm\big(|D|^{-\frac{q}{2}}a|D|^{-\frac{q}{2}}\big) =\tau\left[\big(a_\pm\big)^{\frac{p}{q}}\right]^{\frac{q}{p}},
\end{equation}
where $a_\pm=\frac12(|a|\pm a)$ are the positive and negative parts of $a$ (see~\cite[Proposition~A]{Po:Part1}). 

As mentioned in Section~\ref{sec:Part1.SC}, by the Birman--Schwinger principle these Weyl laws imply the semiclassical Weyl laws~(\ref{eq:Part1.SC-Weyl}). They also lead to the following integration formulas.

\begin{theorem}[{\cite[Theorem~C]{Po:Part1}}]\label{thm:Part1.Integration}
 Assume that Condition~\textup{(W)} is satisfied. For any $a\in  \overline{\sA}$, the operators
$a|D|^{-p}$ and $|D|^{-p/2}a|D|^{-p/2}$, as well as their absolute values, are spectrally measurable, and we have
 \begin{gather}\label{eq:Part1.Int-form1}
 \bint a|D|^{-p} = \bint |D|^{-\frac{p}{2}}a|D|^{-\frac{p}{2}}=\tau\big[a\big],\\
  \bint \big| a|D|^{-p}\big| = \bint \big| |D|^{-\frac{p}{2}}a|D|^{-\frac{p}{2}}\big|=  \tau\left[ |a|\right].
  \label{eq:Part1.Int-form2}
\end{gather}
\end{theorem}

\section{Integration Formulas on Open Manifolds}\label{sec:Cusp.Integration}
In this section, we derive integration formulas for Dirac operators, magnetic Laplacians, and Hodge Laplacians on the open manifold $(X,g)$. 

\subsection{Dirac operators} Suppose that the assumptions of Section~\ref{sec:Cusp.Dirac} hold, i.e., $\overline{X}$ is spin, $g_0$ is an exact cusp metric, and there are no non-zero harmonic spinors on $M$. As shown in Section~\ref{sec:Cusp.Dirac}, the spectral triple~(\ref{eq:Cusp.Dirac-ST}) satisfies Condition~(W) with 
\begin{equation*}
 \tau(f)=\slashed{c}(n)\int_X f(y) d\nu_g(y), \qquad f\in \dot{C}(X)
\end{equation*}
(see~(\ref{eq:Cusp.tau-Dirac})). Therefore, by using Theorem~\ref{thm:Part1.Integration} we get the following integration formulas. 

\begin{theorem}\label{thm:cusp.Int-formulas}
 For every $f\in  \dot{C}(X)$, the operators $f|\shD_g|^{-n}$ and $|\shD_g|^{-n/2}f|\shD_g|^{-n/2}$, as well as their absolute values, are spectrally measurable, and we have
  \begin{gather*}
 \bint f|\shD_g|^{-n} = \bint |\shD_g|^{-\frac{n}{2}}f|\shD_g|^{-\frac{n}{2}}=\slashed{c}(n)\int_X f(y) d\nu_g(y),\\
 \bint \big| f|\shD_g|^{-n}\big| = \bint \big| |\shD_g|^{-\frac{n}{2}}f|\shD_g|^{-\frac{n}{2}}\big|= \slashed{c}(n)\int_X |f(y)| d\nu_g(y). 
   \end{gather*}
\end{theorem}

\subsection{Magnetic and Hodge Laplacians} In what follows, we let $\Delta$ be one of the following operators: 
\begin{itemize}
 \item The magnetic Laplacian $\Delta_A$, where the $1$-form $A$ is trapping (see Section~\ref{sec:Cusp.Magnetic}). In this case we set $c_\Delta:=c(n)$. 
 
 \item The Hodge Laplacian $\Delta_g^{(k)}$ on $k$-forms, where $b_{k-1}(M)=b_k(M)=0$ (see Section~\ref{sec:Cusp.Hodge}). In this case we set $c_\Delta:=c(n,k)$. 
\end{itemize}
As shown in Section~\ref{sec:Cusp.Magnetic} and Section~\ref{sec:Cusp.Hodge}, the spectral triple~(\ref{eq:Cusp.Magnetic-ST}) (resp., (\ref{eq:Cusp.Hodge-ST})) satisfies Condition~(W) with 
\begin{equation*}
 \tau(f)=c_\Delta\int_X f(y) d\nu_g(y), \qquad f\in \dot{C}(X)
\end{equation*}
(see~(\ref{eq:Cusp.tau-Magnetic}) and~(\ref{eq:Cusp.tau-Hodge})). Therefore, by using Theorem~\ref{thm:Part1.Integration} we get the following integration formulas. 

\begin{theorem}\label{thm:Cusp.Integration-Laplacians}
 For every function $f\in  \dot{C}(X)$, the operators $f\Delta^{-n/2}$ and $\Delta^{-n/4}f\Delta^{-n/4}$, as well as their absolute values, are spectrally measurable, and we have
  \begin{gather*}
 \bint f\Delta^{-\frac{n}{2}} = \bint  \Delta^{-\frac{n}{4}} f \Delta^{-\frac{n}{4}} =c_\Delta\int_X f(y) d\nu_g(y),\\
 \bint \big| f\Delta^{-\frac{n}{2}}\big| = \bint \big| \Delta^{-\frac{n}{4}} f \Delta^{-\frac{n}{4}}\big|= c_\Delta\int_X |f(y)| d\nu_g(y). 
   \end{gather*}
\end{theorem}

\begin{remark}
 We refer to~\cite{LU:JMAA25} for strong measurability properties and NC integral formulas for the inverses of a class of magnetic Schr\"odinger operators $\Delta_A+V$ on Euclidean space $\R^n$. Note that in~\cite{LU:JMAA25} the measurability is with respect to the weak Lorentz ideal associated with the function $(t+1)^{-1}\log(t+2)$, which falls outside the scope of this paper. 
\end{remark}

\section{Conformally Cusp Metrics of Infinite Volume}\label{sec:Cusp.infinite}
In this section, we look at conformally cusp metrics $g=x^{2r}g_0$ with $0<r\leq 1/n$. In this case $(X,g)$ is complete and has infinite volume. In what follows, we let $\Delta$ be one of the following operators: 
\begin{itemize}
 \item The square $\shD_g^2$ of the Dirac operator, where $\overline{X}$ is spin, $g_0$ is an exact cusp metric, and there are no non-zero harmonic spinors on $M$ (see Section~\ref{sec:Cusp.Dirac}).
 
 \item The magnetic Laplacian $\Delta_A$, where the $1$-form $A$ is trapping (see Section~\ref{sec:Cusp.Magnetic}). 
 
 \item The Hodge Laplacian $\Delta_g^{(k)}$ on $k$-forms, where $b_{k-1}(M)=b_k(M)=0$ (see Section~\ref{sec:Cusp.Hodge}). 
\end{itemize}
We denote by $\sH$ the corresponding Hilbert space, i.e., $L^2_g(X,\shS_X)$, $L^2_g(X)$, or $L^2_g(X,\Lambda^kT^*X)$. 

The results of~\cite{Mo:MA08, GM:AHP08, GM:TAMS12} that are used in Sections~\ref{sec:Cusp.Dirac}--\ref{sec:Cusp.Hodge} hold for all $r>0$, except for the precise form of the Weyl laws. In particular, the operator $\Delta$ is essentially selfadjoint on $\sH$ and has compact resolvent, and (up to the conjugation by the unitary operator~(\ref{eq:Cusp.isometry}) in the Dirac case) it is a fully elliptic cusp differential operator of type $(2,2r)$. Moreover, in the same way as in Sections~\ref{sec:Cusp.Dirac}--\ref{sec:Cusp.Hodge}, we get a spectral triple, 
\begin{equation}\label{eq:Cusp.infinite-ST}
 \big(\dot{C}^\infty(X), \sH, \sqrt{\Delta}\big). 
\end{equation}
(In the Dirac case we may replace $\sqrt{\Delta}=|\shD_g|$ by $\shD_g$.)

\subsection{The case $r<1/n$} 
Suppose that $r<1/n$. Then, by~\cite[Theorem~3]{Mo:MA08}, \cite[Theorem~4.2]{GM:AHP08}, and \cite[Theorem~1.2]{GM:TAMS12}, there is a constant $\gamma_\Delta>0$ such that, as $j\rightarrow \infty$, we have the Weyl law, 
\begin{equation}\label{eq:Cusp.infinite-Weyl}
 \lambda_j(\Delta) \sim \bigg(\frac{j}{\gamma_\Delta}\bigg)^{2r}. 
\end{equation}
In particular, the spectral triple~(\ref{eq:Cusp.infinite-ST}) is $r^{-1}$-summable. In the Dirac case we have
\begin{equation}\label{eq:Cusp.infinite-gamma}
 \gamma_\Delta= \epsilon_n \frac{\Gamma\big(\frac{1-r}{2r}\big)}{2\sqrt{\pi}\,\Gamma\big(\frac{1}{2r}\big)}\, \zeta\big(\shD_{h_0}; r^{-1}-1\big), 
\end{equation}
where $\zeta(\shD_{h_0}; z):=\Tr[|\shD_{h_0}|^{-z}]$, $\Re z>n-1$, is the zeta function of the Dirac operator of $(M,h_0)$, and $\epsilon_n=1$ if $n$ is odd and $\epsilon_n=2$ if $n$ is even (see~\cite[Theorem~3]{Mo:MA08}). There are similar formulas for magnetic Laplacians and Hodge Laplacians when $g_0$ is an exact cusp metric (see~\cite[Theorem~4.2]{GM:AHP08} and \cite[Eq.~(3.13)]{GM:TAMS12}). Thus, unlike in the case $r>1/n$, the leading coefficient in the Weyl law~(\ref{eq:Cusp.infinite-Weyl}) only depends on the geometry at infinity. 

\begin{lemma}\label{lem:Cusp.infinite-Weyl}
 Assume that $r<1/n$. Let $f\in \dot{C}^\infty(X)$ be~$>0$ on $X$ and have a positive boundary value $f_0$. Then, as $j\rightarrow \infty$, we have
 \begin{equation*}
 \lambda_j\big(f\Delta f\big) \sim \bigg(\frac{j}{\gamma_\Delta f_0^{-1/r}}\bigg)^{2r}. 
\end{equation*}
\end{lemma}
\begin{proof}
 As shown in Sections~\ref{sec:Cusp.Dirac}--\ref{sec:Cusp.Hodge}, the operator $f\Delta f$ is (up to the conjugation by the unitary operator~(\ref{eq:Cusp.isometry}) in the Dirac case) a fully elliptic cusp differential operator of type $(2,2r)$, and its normal operator is given by
 \begin{equation}\label{eq:Cusp.infinite-normal}
 \sN\big(x^{2r}f\Delta f\big)(\xi)= f_0^2\sN\big(x^{2r}\Delta\big)(\xi), \qquad \xi \in \R. 
\end{equation}
Note that this does not require $r$ to be~$>1/n$. Therefore, it follows from~\cite[Propositions~2 and~3]{Mo:MA08} that the function $z\rightarrow \Tr[(f\Delta f)^{-z/2}]$ is holomorphic for $\Re z>r^{-1}$ and has a meromorphic continuation to $\C$. As $r<1/n$, its first pole is located at $z=r^{-1}$. This pole is simple, and the residue is a constant multiple (depending only on $n$ and $r$) of 
\begin{equation*}
 \int_\R \Tr\Big[\sN\big(x^{-1}(f\Delta f)^{-\frac{1}{2r}}\big)(\xi)\Big]d\xi. 
\end{equation*}
As the normal operator map is multiplicative and is compatible with complex powers, it follows from~(\ref{eq:Cusp.infinite-normal}) that $\sN(x^{-1}(f\Delta f)^{-1/2r})(\xi)=f_0^{-1/r}\sN(x^{-1}\Delta^{-1/2r})(\xi)$. Thus, 
\begin{equation*}
 \Res_{z=r^{-1}} \Tr\Big[(f\Delta f)^{-\frac{z}{2}}\Big] = f_0^{-\frac{1}{r}} \Res_{z=r^{-1}} \Tr\Big[\Delta^{-\frac{z}{2}}\Big]. 
\end{equation*}
Combining this with the Tauberian theorem of Delange (see~\cite[Lemma~5]{Mo:MA08}) and the Weyl law~(\ref{eq:Cusp.infinite-Weyl}) gives the result. 
\end{proof}

Lemma~\ref{lem:Cusp.infinite-Weyl} shows that, for $r<1/n$, the spectral triple~(\ref{eq:Cusp.infinite-ST}) satisfies Condition~(W) with $p=r^{-1}$ and 
\begin{equation}\label{eq:Cusp.infinite-tau}
 \tau(f)=\gamma_\Delta f_0, \qquad f\in \dot{C}(X). 
\end{equation}
That is, up to the constant $\gamma_\Delta$, the functional $\tau$ is the evaluation at infinity. Note that $f\rightarrow f_0$ is a character of the $C^*$-algebra $\dot{C}(X)$. In particular, $\tau[\varphi(f)]=\gamma_\Delta\varphi(f_0)$ for all real-valued functions $f\in \dot{C}(X)$ and $\varphi\in C(\R)$. Therefore, by using Theorem~\ref{thm:Part1.SC-Weyl} and Theorem~\ref{thm:Part1.Integration} we arrive at the following result. 

\begin{theorem}\label{thm:Cusp.infinite}
 Assume that $0<r<1/n$, and let $\Delta$ and $\gamma_\Delta$ be as above. 
 \begin{enumerate}
 \item Let $q>0$. For any real-valued potential $V\in \dot{C}(X)$ and every energy level $\lambda \in \R$, we have 
 \begin{equation*}
 \lim_{h\rightarrow 0^+} h^{\frac{1}{r}} N\left(h^{2q}\Delta^q+V;\lambda\right) = \gamma_\Delta \left(V_0-\lambda\right)_{-}^{\frac{1}{2qr}},
\end{equation*}
where $V_0$ is the boundary value of $V$. 

\item For every function $f\in  \dot{C}(X)$, the operators $f\Delta^{-1/2r}$ and $\Delta^{-1/4r}f\Delta^{-1/4r}$, as well as their absolute values, are spectrally measurable, and we have
  \begin{gather*}
 \bint f\Delta^{-\frac{1}{2r}} = \bint  \Delta^{-\frac{1}{4r}} f \Delta^{-\frac{1}{4r}} =\gamma_\Delta f_0,\\
 \bint \big| f\Delta^{-\frac{1}{2r}}\big| = \bint \big| \Delta^{-\frac{1}{4r}} f \Delta^{-\frac{1}{4r}}\big|= \gamma_\Delta |f_0|. 
   \end{gather*}
\end{enumerate}
\end{theorem}

\begin{remark}
 Theorem~\ref{thm:Cusp.infinite} shows that, for $r<1/n$, the semiclassical Weyl laws and the integration formulas only involve the values at infinity. In particular, the semiclassical limit in the first part vanishes when $V_0\geq \lambda$, and in the second part the NC integrals vanish when $f\in C_0(X)$. 
\end{remark}

\subsection{The case $r=1/n$}\label{subsec:Cusp.critical}
Suppose now that $r=1/n$. In this case the results of~\cite{Mo:MA08, GM:AHP08, GM:TAMS12} provide us with \emph{nonclassical} Weyl laws. For instance, in the Dirac case, as $\lambda \rightarrow \infty$, we have
\begin{equation*}
 N\big(|\shD_g|;\lambda\big) \sim  2^{\left[\frac{n}{2}\right]}(2\pi)^{-n}\big|\bS^{n-1}\big| \Vol_{h_0}(M) \lambda^n\log \lambda. 
\end{equation*}
There are similar Weyl laws for magnetic Laplacians and Hodge Laplacians. Therefore, the spectral triple~(\ref{eq:Cusp.infinite-ST}) is not $n$-summable and it does not satisfy Condition~(W), and so the results of~\cite{Po:Part1} do not apply. 

In this case the inverse $|\shD_g|^{-1}$ is contained in a weak Lorentz ideal. It is shown in~\cite{PT:Part1} that $|\shD_g|^{-n}$ is spectrally measurable with respect to the weak Lorentz ideal associated with the function $(t+1)^{-1}\log(t+2)$, and its NC integral with respect to this ideal is equal to
\begin{equation*}
 \frac{1}{n}2^{\left[\frac{n}{2}\right]}(2\pi)^{-n}\big|\bS^{n-1}\big| \Vol_{h_0}(M). 
\end{equation*}
Thus, it recaptures the volume of the boundary $M$. Obtaining semiclassical Weyl laws in the case $r=1/n$ would require extending the results of~\cite{Po:Part1} to spectral triples $(\sA,\sH,D)$ for which $D^{-1}$ lies in a weak Lorentz ideal. We hope to address this question elsewhere.

\end{document}